\documentclass[a4paper,twoside]{article}
\usepackage{amsmath,amsfonts,amsthm,amssymb,amscd}

\makeatletter
\renewenvironment{thebibliography}[1]
     {\section*{\refname
        \@mkboth{\MakeUppercase\refname}{\MakeUppercase\refname}}\small%
      \list{\@biblabel{\@arabic\c@enumiv}}%
           {\settowidth\labelwidth{\@biblabel{#1}}%
            \leftmargin\labelwidth
            \advance\leftmargin\labelsep
            \@openbib@code
            \usecounter{enumiv}%
            \let\p@enumiv\@empty
            \renewcommand\theenumiv{\@arabic\c@enumiv}}%
      \sloppy\clubpenalty4000\widowpenalty4000%
      \sfcode`\.\@m}
     {\def\@noitemerr
       {\@latex@warning{Empty `thebibliography' environment}}%
      \endlist}
\newcommand{\leftsideset}[2]{%
  \@mathmeasure\z@\displaystyle{#2}%
  \global\setbox\@ne\vbox to\ht\z@{}\dp\@ne\dp\z@
  \setbox\tw@\box\@ne
  \@mathmeasure4\displaystyle{\copy\tw@#1}%
  \@mathmeasure6\displaystyle{#2}%
  \dimen@-\wd6 \advance\dimen@\wd4 \advance\dimen@\wd\z@
  \hbox to\dimen@{}\mathord{\kern-\dimen@\box4\box6}%
}
\makeatother

\let\ii=\i
\catcode`\ˆ=\active\defˆ{\`a}
\catcode`\=\active\def{\`e}
\catcode`\Ž=\active\defŽ{\'e}
\catcode`\'=\active\def'{\'{\ii}}
\catcode`\˜=\active\def˜{\`o}
\catcode`\œ=\active\defœ{\'u}
\catcode`\é=\active\defé{\`E}
\catcode`\£=\active
\newcommand{£}{\pounds\ida}
\newcommand{\Lie}{{\,\,\tilde{\!\!\pounds}}\ida}
\DeclareMathAlphabet{\mathsf}{OT1}{cmss}{m}{n}

\DeclareMathOperator{\im}{im}

\DeclareMathOperator{\tr}{tr}
\newcommand{\thorn}{\hbox{\textsf{\setbox0=\hbox{l}\copy0\kern-\wd0 p}}}

\newcommand{\Aut}{\mathrm{Aut}}

\newcommand{\Ad}{\mathrm{Ad}}

\renewcommand{\ast}{{\displaystyle*}}

\newcommand{\C}{\mathbb{C}}

\newcommand{\cde}{\nabla\ida}
\newcommand{\cf}{\protect\emph{cf.{}}}

\renewcommand{\d}{\mathrm{d}}
\newcommand{\de}{\partial}
\def\dual(#1,#2){\langle{#1},{#2}\rangle}

\newcommand{\g}{\mathfrak{g}}
\newcommand{\ga}{\gamma}

\newcommand{\Gadash}{$\Gamma$\nobreakdash-\hspace{0pt}}

\newcommand{\Gdash}{$G$\nobreakdash-\hspace{0pt}}
\newcommand{\gl}{\mathfrak{gl}}
\newcommand{\GL}{\mathrm{GL}}
\newcommand{\Gkdash}{$G^k_m$\nobreakdash-\hspace{0pt}}
\newcommand{\GLdash}{$\GL(m,\R)$\nobreakdash-\hspace{0pt}}
\newcommand{\h}{\mathfrak{h}}
\newcommand{\Hdash}{$H$\nobreakdash-\hspace{0pt}}

\renewcommand{\i}{\mathrm{i}}
\newcommand{\id}{\mathrm{id}}
\newcommand*{\ida}[2]{\ifx#1^{}^{#2}
                      \else\ifx#1_{}_{\!#2}
                      \else\errmessage{Sub/Superscript token missing}\fi
                      \fi}
\newcommand{\KK}{\mathsf{K}}
\renewcommand{\l}{\lambda}

\newcommand{\m}{\mathfrak{m}}
\newcommand{\mdash}{$m$\nobreakdash-\hspace{0pt}}

\newcommand{\onedash}{$1$\nobreakdash-\hspace{0pt}}

\newcommand{\pr}{\mathrm{pr}}
\newcommand{\R}{\mathbb{R}}

\newcommand{\Symm}{{\textstyle\bigvee}}

\renewcommand{\sl}{\mathfrak{sl}}
\newcommand{\SL}{\mathrm{SL}}
\newcommand{\so}{\mathfrak{so}}
\newcommand{\SO}{\mathrm{SO}}
\newcommand{\SOdash}{$\SO(p,q)$\nobreakdash-\hspace{0pt}}

\newcommand{\Spin}{\mathrm{Spin}}
\newcommand{\Spindash}{$\Spin(p,q)$\nobreakdash-\hspace{0pt}}

\newcommand{\Tw}{\leftsideset{^w}T}

\newcommand{\V}{\mathbb{V}}

\newcommand{\WGdash}{$W^{k,h}_mG$\nobreakdash-\hspace{0pt}}
\newcommand{\X}{\mathfrak{X}}

\newcommand{\xiK}{\xi_{\mathrm{K}}{}}

\renewcommand{\[}{\left[}
\renewcommand{\]}{\right]}
\newtheorem{proposition}{Proposition}[section]
\newtheorem{corollary}[proposition]{Corollary}
\newtheorem{lemma}[proposition]{Lemma}
\newtheorem{theorem}[proposition]{Theorem}
\theoremstyle{definition}
\newtheorem{definition}[proposition]{Definition}
\newtheorem{example}[proposition]{Example}
\newtheorem{remark}[proposition]{Remark}
\numberwithin{equation}{section}

\begin{document}

\hyphenation{ac-know-ledges South-amp-ton there-by}

\title{Reductive \Gdash structures and Lie derivatives}
\date{22$^{\text{nd}}$ August 2003}

\author{Marco Godina\thanks{Dipartimento di Matematica, Universitˆ di Torino, Via Carlo Alberto 10, 10123
Torino, Italy. E-mail address: godina@matlag.dm.unito.it.}\ \ and Paolo Matteucci\thanks{Faculty of
Mathematical Studies, University of Southampton, Highfield, Southampton SO17~1BJ, UK. 
E-mail address: p.matteucci@maths.soton.ac.uk.}}

\maketitle
\begin{abstract}
Reductive \Gdash structures on a principal bundle $Q$ are considered. It is shown that these structures,
i.e.\ reductive \Gdash subbundles $P$ of $Q$, admit a canonical decomposition of the pull-back vector bundle
$i_{P}^{\ast}(TQ)\equiv P\times_{Q}TQ$ over~$P$. For classical \Gdash structures, i.e.\ reductive \Gdash
subbundles of the linear frame bundle, such a decomposition defines an infinitesimal canonical lift. This
lift extends to a prolongation \Gadash structure on $P$. In this general geometric framework the theory of
Lie derivatives is considered. Particular emphasis is given to the morphisms which must be taken in order to
state what kind of Lie derivative has to be chosen. On specializing the general theory of gauge-natural Lie
derivatives of spinor fields to the case of the Kosmann lift, we recover the result originally found by
Kosmann. We also show that in the case of a reductive \Gdash structure one can introduce a ``reductive Lie 
derivative'' with respect to a certain class of generalized infinitesimal automorphisms. This differs, in
general, from the gauge-natural one, and we conclude by showing that the ``metric Lie derivative'' introduced
by Bourguignon and Gauduchon is in fact a particular kind of reductive rather than gauge-natural Lie
derivative.
\end{abstract}

\pagestyle{myheadings}
\markboth{Marco Godina and Paolo Matteucci}%
         {Reductive \protect\Gdash structures and Lie derivatives}
\thispagestyle{plain}

\section*{Introduction}
It has now become apparent that there has been some confusion regarding the concept of a Lie derivative of
spinor fields, both in the mathematical and the physical literature.

Lichnerowicz was the first one to give a correct definition for such an object, although with respect to
infinitesimal isometries only. The local expression given by Lichnerowicz in 1963 \cite{lichner63} is
\begin{equation}
£_{\xi}\psi := \xi^{a}\cde_{a}\psi
-\frac14\nabla_{a}\xi_{b}
\gamma^{a}\gamma^{b}\psi,
\tag{$*$}
\label{eq:lld}
\end{equation}
where $\cde_{a}\xi_{b}=\cde_{[a}\xi_{b]}$,
as $\xi$ is assumed to be a Killing vector field. 

After a first attempt to extend Lichnerowicz's definition to generic infinitesimal
transformations \cite{kosmann66}, in 1972 Kosmann put forward a new definition of a Lie
derivative of spinor fields in \cite{kosmann72}, her doctoral thesis under Lichnerowicz's supervision.
Indeed, in her previous work she had just extended \emph{tout court} Lichnerowicz's definition to the case
of a generic vector field~$\xi$, without antisymmetrizing $\cde_{a}\xi_{b}$. Therefore, the local expression
appearing in \cite{kosmann66} could not be given any clear-cut geometrical meaning. The remedy was then
realized to be retaining Lichnerowicz's local expression~\eqref{eq:lld} for a \emph{generic} vector
field~$\xi$, but explicitly taking the antisymmetric part of $\cde_{a}\xi_{b}$ only \cite{kosmann72}.

Several papers on the subject followed, including particularly Binz and Pferschy's \cite{bp72} and
Bourguignon and Gauduchon's \cite{bg92}. Furthermore, among the physics community much interest has been
attracted by Penrose and Rindler's definition \cite{pr2}, despite its being restricted to infinitesimal
conformal isometries because of the (implicit) requirement that the Lie derivative commute with the
isomorphism between the complexified tangent bundle and the tensor product of the spinor bundle and its
complex conjugate (see \cite{delaney} for a thorough discussion).

In this paper we investigate whether the definition of a Lie derivative of spinor
fields can be placed in the more general framework of the theory of Lie derivatives of sections of
fibred manifolds (and, more generally, of differentiable maps between two manifolds) stemming from
Trautman's 1972 seminal paper \cite{trautman72} and further developed by Jany\v ska and
Kol\'a\v{r}~\cite{jk82} (see also \cite{kms}).

A first step in this direction was already taken in \cite{fffg96}, where Kosmann's 1972 definition was
successfully placed in the framework of the theory of Lie derivatives of sections of \emph{gauge-natural
bundles} by introducing a new geometric concept, which the authors called the ``Kosmann lift''.

The aim of this paper is to provide a more transparent geometric explanation of the Kosmann lift and, at
the same time, a generalization to reductive \Gdash structures. Indeed, the Kosmann lift is but a
\emph{particular case} of this interesting generalization.

The structure of the paper is as follows: in \S\ref{sec:notation} preliminary notions on
principal bundles are recalled for the main purpose of fixing our notation; in \S\ref{sec:rGs}
the concept of a reductive \Gdash structure and its main properties are introduced; in
\S\ref{sec:gnb} a constructive approach to gauge-natural bundles is proposed together with a
number of relevant examples; in \S\ref{sec:ss} split structures on principal bundles are
considered and the notion of a generalized Kosmann lift is defined; finally, in \S\ref{sec:ld}
the general theory of Lie derivatives is applied to the context of reductive \Gdash structures,
allowing us to analyse the concept of the Lie derivative of spinor fields in all its different
flavours from the most general point of view. The proofs of the results presented in this paper
mainly consist of the careful application of the definitions which precede them, and therefore
are mostly omitted.

\section{Notation}\label{sec:notation}
Let $M$  be a manifold and $G$ a Lie group. A \emph{principal} (\emph{fibre}) \emph{bundle}~$P$  over~$M$
with structure group~$G$  is obtained by attaching a copy of~$G$ to each point of~$M$, i.e.\ by giving a
\Gdash manifold~$P$,  on which~$G$ acts on the right and which satisfies the following conditions:
\begin{enumerate}
\item The (right) action $r\colon P\times G\to P$ of~$G$ on~$P$ is \emph{free}, i.e.\ $u\cdot a :=r(u,a) =u$, $u\in P$, implies $a=e$, $e$ being the unit element of~$G$.

\item $M=P/G$ is the quotient space of~$P$ by the equivalence relation induced by~$G$, i.e.~$M$ is the space
of orbits. Moreover, the canonical projection $\pi\colon P\rightarrow M$ is smooth.

\item $P$ is locally trivial, i.e.~$P$ is locally a product $U\times G$,
where~$U$ is an open set in~$M$. More precisely, there exists a diffeomorphism
$\Phi\colon\pi^{-1}(U)\to U\times G$ such that $\Phi(u)=(\pi(u), f(u))$,
where the mapping $ f\colon\pi^{-1}(U)\to G$ is \emph{\Gdash equivariant},
i.e.\ $ f(u\cdot a)= f(u)\cdot a$ for all $u \in \pi^{-1}(U)$,
$a \in G$.
\end{enumerate}
A principal bundle will be denoted by $(P,M,\pi;G)$, $P(M,G)$,
$\pi\colon P\to M$ or simply $P$, according to the particular context.
$P$ is called the \emph{bundle} (or \emph{total}) \emph{space}, $M$ the \emph{base}, $G$ the
\emph{structure group}, and $\pi$ the \emph{projection}. The closed submanifold $\pi^{-1}(x)$, $x\in
M$, will be called the \emph{fibre over}~$x$. For any point $u\in P$,
we have $\pi^{-1}(x)=u\cdot G$, where $\pi(u)=x$, and $u\cdot G$ will be called
the \emph{fibre through}~$u$. Every fibre is diffeomorphic to~$G$,
but such a diffeomorphism depends on the chosen trivialization.

Given a manifold~$M$ and a Lie group~$G$, the product manifold $M\times G$ is
a principal bundle over~$M$ with projection $\pr_{1}\colon M\times G\to M$
and structure group~$G$, the action being given by $(x,a)\cdot b= (x,a\cdot b)$.
The manifold $M\times G$ is called a \emph{trivial principal bundle}.

A \emph{homomorphism} of a principal bundle $P'(M',G')$ into another
principal bundle \allowbreak $P(M,G)$ consists of a differentiable mapping
$\Phi\colon P'\to P$ and a Lie group homomorphism $ f\colon G'\to G$
such that $\Phi(u'\cdot a')=\Phi(u')\cdot  f(a')$
for all $u' \in P'$, $a' \in G'$. Hence, $\Phi$  maps fibres into fibres
and induces a differentiable mapping $\varphi\colon M'\to M$ by
$\varphi(x')=\pi(\Phi(u'))$, $u'$ being an arbitrary point over~$x'$.
A homomorphism $\Phi\colon P'\to P$ is called an \emph{embedding} if
$\varphi\colon M'\to M$ is an embedding and $ f\colon G'\to G$
is injective. In such a case, we can identify~$P'$ with~$\Phi(P')$,
$G'$ with~$ f(G')$ and~$M'$ with~$\varphi(M')$, and~$P'$ is said to be
a \emph{subbundle} of~$P$. If $M'=M$ and $\varphi=\id_M$, $P'$ is called
a \emph{reduced subbundle} or a \emph{reduction} of~$P$, and we also say that~$G$
``reduces'' to the subgroup~$G'$.

A homomorphism $\Phi\colon P'\to P$ is called an \emph{isomorphism} if there
exists a homomorphism of principal bundles $\Psi\colon P\to P'$ such that
$\Psi \circ \Phi=\id_{P'}$ and $\Phi \circ \Psi=\id_{P}$.

\section{Reductive $G$-structures and their prolongations}\label{sec:rGs}
\begin{definition}
Let $H$ be a Lie group and $G$ a Lie subgroup of $H$\textup.
Denote by~$\h$ the Lie algebra of~$H$ and by~$\g$
the Lie algebra of~$G$\textup. We shall say that $G$ is a \emph{reductive 
Lie subgroup} of~$H$ if there exists a direct sum decomposition
\begin{equation*}
\h = \g \oplus \m,
\label{eq:GA-split}
\end{equation*}
where~$\m$ is an $\Ad_G$-invariant vector subspace of~$\h$\textup,
i\textup.e\textup.\ $\Ad_a(\m) \subset \m$ for all $a\in G$ 
\textup(which means that the $\Ad_G$ representation of $G$ in $\h$ is reducible
into a direct sum decomposition of two $\Ad_G$-invariant vector spaces\textup:
\cf~\textup{\cite{kn1},} p\textup.~$83$\textup{).}
\end{definition}
\begin{remark}
A Lie algebra $\h$ and a Lie subalgebra $\g$
satisfying these properties form a so-called \emph{reductive pair}
(cf. \cite{cb2}, p.\ 103). Moreover, $\Ad_G(\m) \subset \m$ implies
$\[\g , \m \] \subset \m$, and, conversely, if $G$ is connected,
$\[\g , \m \] \subset \m$ implies $\Ad_G(\m) \subset \m$ 
\cite[p.\ 190]{kn2}.
\end{remark}
\begin{example}
Consider a subgroup $G \subset H$ and suppose
that an $\Ad_{G}$-invariant metric $K$ can be assigned on the
Lie algebra~$\h$ (e.g., if $H$ is a semisimple Lie group, $K$
could be the Cartan-Killing form: indeed, this form is $\Ad_{H}$-invariant
and, in particular, also $\Ad_{G}$-invariant). Set
\begin{equation*}
\m :=\g^{\perp} \equiv\{\, v\in\h \mid K(v,u) = 0\ \forall u\in\g \,\} \, .
\end{equation*}
Obviously, $\h$ can be decomposed as the direct sum $\h = \g \oplus \m$
and it is easy to show that, under the assumption of $\Ad_{G}$-invariance
of $K$, the vector subspace $\m$ is also $\Ad_{G}$-invariant.
\end{example}
\begin{example}[The unimodular group]\label{ex:SL}
The unimodular group $\SL(m,\R)$ is an example of a reductive 
Lie subgroup of $\GL(m,\R)$. To see this, first recall
that its Lie algebra $\sl(m,\R)$ is formed by all $m\times m$ traceless matrices.
If $\mathsf M$ is any matrix in $\gl(m,\R)$, the following decomposition holds:
\begin{equation*}
\mathsf M =\mathsf U +1/m\tr(\mathsf M)\mathsf I,
\end{equation*}
where $\mathsf I := \id_{\gl(m, \R)}$ and $\mathsf U$ is traceless. Indeed,
\begin{equation*}
\tr(\mathsf U) = \tr(\mathsf M)  -1/m\tr(\mathsf M)\tr(\mathsf I) = 0.
\end{equation*}
Accordingly, the Lie algebra $\gl(m,\R)$ can be decomposed as follows:
\begin{equation*}
\gl(m, \R) =\sl(m,\R) \oplus\R\mathsf I.
\end{equation*}
In this case, $\m$ is the set of all real multiples of~$\mathsf I$, 
which is obviously adjoint-invariant under $\SL(m, \R)$.
Indeed, if $\mathsf S$ is an arbitrary element of $\SL(m, \R)$, 
for any $a\in\R$ one has
\begin{equation*}
\Ad_{\mathsf S}(a\mathsf I) \equiv\mathsf S(a\mathsf I)\mathsf S^{-1} 
                            =a\mathsf{ISS}^{-1} = a\mathsf I.
\end{equation*}
This proves that $\R\mathsf I$ is adjoint-invariant under $\SL(m, \R)$,
and $\SL(m,\R)$ is a reductive Lie subgroup of $\GL(m, \R)$.
\end{example}
Given the importance of the following example for the future developments of the theory, we shall
state it as
\begin{proposition}The \textup(pseudo\nobreakdash-\textup) orthogonal group $\SO(p,q)$\textup,
$p+q=m$\textup, is  a reductive Lie subgroup of $\GL(m,\R)$\textup.
\end{proposition}
\begin{proof} 
Let $\eta$ denote the standard metric of signature
$(p,q)$, with $p+q=m$, on $\R^{m}\equiv\R^{p,q}$ and $\mathsf M$ be any matrix in $\gl(m,\R)$.
Denote by $\mathsf M^\top$ the adjoint (``transpose'') of~$\mathsf M$ with respect to~$\eta$,
defined by requiring $\eta(\mathsf M^\top v,v')=\eta(v,\mathsf Mv')$ for all $v, v'\in\R^m$. Of
course, any traceless matrix can be (uniquely) written as the sum of an antisymmetric matrix and
a symmetric traceless matrix. Therefore,
\begin{equation*}
\sl(m, \R) =\so(p,q) \oplus\V,
\end{equation*}
$\so(p,q)$ denoting the Lie algebra of the (pseudo-) orthogonal group $\SO(p,q)$ for~$\eta$, formed by
all matrices~$\mathsf A$ in $\gl(m,\R)$ such that $\mathsf A^\top=-\mathsf A$, and $\V$ the vector space
of all matrices~$\mathsf V$ in $\sl(m,\R)$ such that $\mathsf V^\top=\mathsf V$. Now, let
$\mathsf O$ be any element of $\SO(p,q)$ and set $\mathsf V':=\Ad_{\mathsf O}\mathsf V\equiv
\mathsf {OVO}^{-1}$ for any $\mathsf V\in\V$. We have
\begin{equation*}
\mathsf V'{}^\top =(\mathsf {OVO}^\top)^\top= \mathsf V'
\end{equation*}
because $\mathsf V^\top=\mathsf V$ and $\mathsf O^{-1}=\mathsf O^\top$. Moreover,
\begin{equation*}
\tr(\mathsf V') =\tr(\mathsf O)\tr(\mathsf V)\tr(\mathsf O^{-1})=0
\end{equation*}
since $\mathsf V$ is traceless. So, $\mathsf V'$ is in~$\V$, thereby proving that $\V$ is
adjoint-invariant under $\SO(p,q)$. Therefore, $\SO(p,q)$ is a reductive Lie subgroup of
$\SL(m,\R)$ and, hence, also a reductive Lie  subgroup of $\GL(m, \R)$ by virtue of
Example~\ref{ex:SL}.
\end{proof}
\begin{definition}
A \emph{reductive \Gdash structure}
on a principal  bundle $Q(M, H)$ is a principal
subbundle $P(M, G)$ of $Q(M, H)$ such that $G$ is a reductive Lie subgroup
of~$H$\textup.
\end{definition}
Now, since later on we shall consider the case of spinor fields, it is convenient to give the following
general
\begin{definition}
Let $P(M,G)$ be a principal bundle and $\rho\colon \Gamma \to G$ a central homomorphism of a Lie 
group $\Gamma$ onto~$G$, i.e.\ such that its kernel is discrete and contained in the centre of~$\Gamma$
\cite{gp78} (see also \cite{hae56}). A \emph{\Gadash structure} on $P(M,G)$ is a principal bundle map
$\zeta \colon \tilde P \to P$ which is equivariant under the right actions of the structure groups\textup,
i\textup.e\textup. 
\begin{equation*}
\zeta(\tilde u\cdot\alpha) =\zeta(\tilde u)\cdot\rho(\alpha)
\end{equation*}
for all $\tilde u\in\tilde P$ and $\alpha\in\Gamma$\textup.
\end{definition}
Equivalently, we have the following commutative diagrams
\begin{equation*}
\begin{CD}
   \tilde P @>\zeta>>  P  \\
   @V{\tilde\pi}VV    @VV{\pi}V  \\
   M        @>>\id_M> M
\end{CD}
\qquad\qquad
\begin{CD}
   \tilde P @>\tilde r^\alpha>> \tilde P  \\
   @V{\zeta}VV                   @VV{\zeta}V  \\
   P        @>>r^a>  P
\end{CD}
\end{equation*}
$r^a$ and $\tilde r^\alpha$ denoting the 
right action on~$P$ and~$\tilde P$, respectively (see \cite{fffg98}).
This means that, for $\tilde u\in\tilde P$, both $\tilde u$ and 
$\zeta(\tilde u)$ lie over the same point, and~$\zeta$, restricted to 
any fibre, is a ``copy'' of $\rho$, i.e.\ it is equivalent to it.
The existence condition for a \Gadash structure on $P$ can be
formulated in terms of {\v C}ech cohomology \cite{hae56,gp78,lm}.
\begin{remark}
The bundle map $\zeta \colon \tilde P \to P$ is a covering space since its kernel is discrete.
\end{remark}
Recall now that for any principal bundle
$(P,M,\pi,G)$ a (\emph{principal}) \emph{automorphism} of $P$ is a diffeomorphism
$\Phi \colon P \to P$ such that $\Phi (u \cdot a) = \Phi (u) \cdot a$ for every $u \in P$, $a \in G$.
Each $\Phi$ induces a unique
diffeomorphism $\varphi\colon M\to M$ 
such that $\pi\circ\Phi =\varphi\circ\pi$. Accordingly, we shall denote by
$\Aut(P)$ the group of all principal automorphisms of $P$.
Assume that a vector field $\Xi$ on $P$ generates
a local \onedash parameter group $\{\Phi_t\}$. Then, $\Xi$ is \emph{\Gdash invariant}
if and only if $\Phi_t$ is an automorphism of $P$ for every $t\in\R$.
Accordingly, we denote by $\X_G(P)$
the Lie algebra of \Gdash invariant vector fields on $P$.

Now recall that, given a fibred manifold $\pi\colon B\to M$, 
a \emph{projectable vector field} on~$B$ over a vector field~$\xi$ on~$M$ is a vector
field~$\Xi$  on~$B$ such that $T\pi\circ\Xi=\xi\circ\pi$. It follows
\begin{proposition}
Let $P(M,G)$ be a principal bundle\textup.
Then\textup, every \Gdash invariant vector field $\Xi$ on $P$ is projectable
over a unique vector field $\xi$ on the base manifold $M$\textup.
\end{proposition}
\begin{proposition}\label{prop:GivfGaivf}
Let $\zeta \colon \tilde P \to P$ be a \Gadash structure on $P(M,G)$\textup. Then\textup, every \Gdash
invariant vector field $\Xi$ on $P$ admits  a unique \textup(\Gadash invariant\textup) lift $\tilde\Xi$
onto~$\tilde P$\textup.
\end{proposition}
\begin{proof} 
Consider a \Gdash invariant vector field $\Xi$, its flow being denoted by
$\{\Phi_{t}\}$.
For each $t\in\R$, $\Phi_t$ is an automorphism of~$P$. Moreover, $\zeta \colon \tilde P \to P$ being
a covering space, it is possible to lift $\Phi_{t}$ to a (unique) bundle map $\tilde\Phi_t\colon\tilde P
\to\tilde P$ in the following way. For any point $\tilde u\in\tilde P$, consider the (unique) point
$\zeta(\tilde u) =u$. From the theory of covering spaces it follows that, for the curve
$\gamma_u\colon\R\to P$ based at~$u$, that is $\gamma_u(0)=u$, and defined by $\gamma_u(t) :=\Phi_t(u)$,
there exists a unique curve $\tilde \gamma_{\tilde u}\colon\R\to\tilde P$ based at~$\tilde u$
such that $\zeta \circ \tilde \gamma_{\tilde u} =\gamma_{u}$. It is possible to define a principal bundle map
$\tilde\Phi_t\colon\tilde P\to\tilde P$ covering $\Phi_t$ by setting $\tilde\Phi_t(\tilde u) :=\tilde
\gamma_{\tilde u}(t)$. The \onedash parameter group of automorphisms $\{\tilde\Phi_t\}$ of $\tilde P$
defines a vector field $\tilde\Xi(\tilde u) :=\frac{\de}{\de t}[\tilde\Phi_t(\tilde u)]\big\rvert_{t=0}$
for all $\tilde u\in\tilde P$.
\end{proof} 
\begin{proposition}
Let $\zeta \colon \tilde P \to P$ be a \Gadash structure on $P(M,G)$\textup. Then\textup, every \Gadash
invariant vector field $\tilde\Xi$ on $\tilde P$  is projectable over a unique \Gdash invariant vector field
$\Xi$ on~$P$\textup.
\end{proposition}
\begin{proof} 
Consider a \Gadash invariant vector field $\tilde\Xi$ on $\tilde P$. Denote its flow by $\{\tilde\Phi_t\}$.
Each $\tilde\Phi_t$ induces a unique automorphism $\Phi_t\colon P\to P$ such that $\zeta \circ \tilde\Phi_t =
\Phi_t\circ\zeta$ and, hence, a unique vector field $\Xi$ on $P$ given by  $\Xi(u)
:=\frac{\de}{\de t}[\Phi_t(u)]\big\rvert_{t = 0}$ for all $u\in P$.
\end{proof} 
\begin{corollary}
Let $\zeta \colon \tilde P \to P$ be a \Gadash structure on $P(M,G)$\textup.
There is a bijection between \Gdash invariant vector fields on $P$
and \Gadash invariant vector fields on $\tilde P$\textup.
\end{corollary}

\section{Gauge-natural bundles}\label{sec:gnb}
In this section we shall introduce the category of gauge-natural bundles \cite{eck81,kms} and give a number
of relevant examples. Geometrically,  gauge-natural bundles possess a very rich structure, which generalizes
the classical one of natural bundles. From the physical point of view, this framework enables one to
treat at the same time, under a unifying formalism, natural field theories such as general relativity, gauge
theories, as well as bosonic and fermionic matter field theories (\cf\ \cite{fffg98,faf98,gmv,gng}).
\begin{definition}
Let $j^\ell_pf$ denote the $\ell$-th order jet prolongation 
of a map~$f$ evaluated at a point~$p$\textup. The set
\begin{equation*}
\{\,j^k_0\alpha \mid \alpha\colon\R^{m} \to \R^{m}\text{\textup, } 
\alpha(0)=0 \text{\textup, locally invertible}\,\}
\end{equation*}
equipped with the jet composition $j^k_0\alpha\circ j^k_0\alpha':=j^k_0(\alpha\circ\alpha')$ is a
Lie group called the \emph{$k$-th differential group} and denoted by~$G^k_m$.
\end{definition}
For $k=1$ we have, of course, the identification $G^1_m\cong\GL(m,\R)$.
\begin{definition}
Let $M$ be an \mdash dimensional manifold\textup. The principal bundle over~$M$ with group $G^k_m$
is called the \emph{$k$-th order frame bundle} over~$M$ and will be denoted by $L^k\!M$\textup.
\end{definition}
For $k=1$ we have, of course, the identification $L^1\!M\cong LM$, where $LM$ is the usual
(principal) \emph{bundle of linear frames} over~$M$ (\cf, e.g., \cite{kn1}).
\begin{definition}
Let $G$ be a Lie group\textup.
Then\textup, the \emph{space of $(m,h)$-velocities} of~$G$ is defined as
\begin{equation*}
T^h_mG:=\{\,j^h_0a\mid a\colon\R^{m} \to G\,\}.
\end{equation*}
\end{definition}
Thus, $T^h_mG$ denotes the set of $h$-jets with source at the origin $0 \in \R^m$ and target in~$G$. It is a
subset of the manifold $J^h(\R^m, G)$ of $r$-jets with source in~$\R^m$ and target in~$G$. The set
$J^h(\R^m, G)$ is a fibre bundle over~$\R^m$ with respect to the canonical jet projection of $J^h(\R^m, G)$
on~$\R^m$, and $T^h_mG$ is its fibre over $0\in\R^m$. Moreover, the set $T^h_mG$ can be given the 
structure of a Lie group. Indeed, let $S,T\in T^h_mG$ be any elements.  We define a (smooth) multiplication
in $T^h_mG$ by:
\begin{equation*}
\left\{
\begin{aligned}
{}&T^h_m\mu\colon T^h_mG \times T^h_mG \to T^h_mG  \\
{}&T^h_m\mu\colon (S=j^h_0a,T=j^h_0b)\mapsto
S\cdot T:= j^{h}_{0}(a\cdot b)
\end{aligned}
\right.,
\end{equation*}
where $(a\cdot b)(x):=a(x)\cdot b(x)\equiv\mu(a(x),b(x))$ 
is the group multiplication in~$G$.
The mapping $(S,T)\mapsto S\cdot T$ is associative;
moreover, the element $j^h_0e$, $e$ denoting both the unit element in~$G$
and the constant mapping from~$\R^m$ to~$e$,
is the unit element of $T^h_mG$,
and $j^r_0a^{-1}$, where $a^{-1}(x):=\big(a(x)\big)^{-1}$ 
(the inversion being taken in the group $G$),
is the inverse of~$j^h_0a$.
\begin{definition}
Consider a principal bundle $P(M,G)$. 
Let $k$ and $h$ be two natural numbers such that $k\geq h$\textup.
Then\textup, by the \emph{$(k,h)$-principal prolongation} of~$P$ we shall mean the bundle
\begin{equation}
W^{k,h}\!P:= L^k\!M\times_M J^h\!P,
\label{eq:WP}
\end{equation}
where $ L^k\!M$ is the $k$-th order frame bundle of~$M$ and $J^h\!P$ denotes the
$h$-th order jet prolongation of~$P$\textup.
A point of $W^{k,h}\!P$ is of the form $(j^k_0\epsilon,j^h_x\sigma)$, where 
$\epsilon\colon\R^{m} \to M$ is locally invertible and such that $\epsilon(0)=x$, and 
$\sigma\colon M\to P$ is a local section around the point $x\in M$\textup.
\end{definition}
Unlike $J^h\!P$, $W^{k,h}\!P$ is a principal
bundle over~$M$ whose structure group is
\begin{equation*}
W^{k,h}_mG:= G^k_m \rtimes T^h_mG.
\end{equation*}
$W^{k,h}_mG$ is called the \emph{$(m;k,h)$-principal prolongation} of~$G$. The group
multiplication on
$W^{k,h}_mG$ is defined by the following rule:
\begin{equation*}
(j^k_0\alpha,j^h_0a)\odot(j^k_0\beta,j^h_0b):=
\Big(j^k_0(\alpha\circ\beta),j^h_0\big((a\circ\beta)\cdot b\big)\Big),
\end{equation*}
`$\cdot$' denoting the group multiplication in~$G$.
The right action of $W^{k,h}_mG$ on $W^{k,h}\!P$
is then defined by:
\begin{equation*}
(j^k_0\epsilon,j^h_x\sigma)\odot(j^k_0\alpha,j^h_0a):=
\Big(j^k_0(\epsilon\circ\alpha),
j^h_x\big(\sigma\cdot(a\circ \alpha^{-1}\circ\epsilon^{-1})\big)\Big),
\end{equation*}
`$\cdot$' denoting now the canonical right action of~$G$ on~$P$.
\begin{definition}
Let $\Phi\colon P\to P$ be an automorphism over a diffeomorphism
$\varphi\colon M\to M$\textup. We define an \emph{automorphism of $W^{k,h}\!P$
associated with~$\Phi$} by
\begin{equation}
\left\{
\begin{aligned}
{}&W^{k,h}\Phi\colon W^{k,h}\!P\to W^{k,h}\!P  \\
{}&W^{k,h}\Phi\colon(j^k_0\epsilon,j^h_x\sigma)\mapsto
\big(j^k_0(\varphi\circ\epsilon),j^h_x(\Phi\circ\sigma\circ\varphi^{-1})\big)
\end{aligned}
\right..
\label{eq:WPhi}
\end{equation}
\end{definition}
\begin{proposition}
The bundle morphism $W^{k,h}\Phi$ preserves the right action\textup, thereby being a principal
automorphism\textup.
\end{proposition}
By virtue of~\eqref{eq:WP} and~\eqref{eq:WPhi} $W^{k,h}$ turns out to be a functor from the category	of
principal \Gdash bundles over \mdash dimensional manifolds and local isomorphisms to the category of
principal \WGdash bundles \cite{kms}. Now, let $P_\lambda:=W^{k,h}\!P\times_{\lambda} F$ be a fibre bundle 
associated with $P(M,G)$ via an action~$\lambda$ of $W^{k,h}_mG$ on a manifold~$F$. There exists canonical
representation of the automorphisms of~$P$ induced by~\eqref{eq:WPhi}. Indeed, if
$\Phi\colon P\to P$ is an automorphism over a diffeomorphism $\varphi\colon M\to M$, then we can define the
corresponding 
\emph{induced automorphism}~$\Phi_\lambda$ as
\begin{equation}
\left\{
\begin{aligned}
{}&\Phi_\l\colon P_\l\to P_\l  \\
{}&\Phi_\l\colon[u,f]_\l\mapsto[W^{k,h}\Phi(u),f]_\l
\end{aligned}
\right.,
\label{eq:indaut}
\end{equation}
which is well-defined since it turns out to be independent of the representative $(u,f)$, $u\in P$,
$f\in F$. This construction yields a functor $\cdot_\l$ from the category of principal \Gdash bundles to
the category of fibred manifolds and fibre-respecting mappings.
\begin{definition}
A \emph{gauge-natural bundle of order~$(k,h)$} over~$M$ associated with $P(M,G)$ is any such functor\textup.
\end{definition}
If we now restrict attention to the case $G=\{e\}$ and $h=0$, we can recover the classical notion of natural
bundles over~$M$. In particular, we have the following
\begin{definition}
Let $\varphi\colon M\to M$ be a diffeomorphism\textup. We define an automorphism of $L^k\!M$
associated with~$\varphi$\textup, called its \emph{natural lift}\textup, by
\begin{equation*}
\left\{
\begin{aligned}
{}&L^k\varphi\colon L^k\!M\to L^k\!M  \\
{}&L^k\varphi\colon j^k_0\epsilon\mapsto j^k_0(\varphi\circ\epsilon)
\end{aligned}
\right..
\end{equation*}
\end{definition}
Then, $L^k$ turns out to be a functor from the category of \mdash dimensional manifolds and local
diffeomorphisms to the category of principal \Gkdash bundles. Now, given any fibre bundle
associated with $L^k\!M$ and any diffeomorphism on~$M$, we can define a corresponding induced automorphism
in the usual fashion. This construction yields a functor from the category of \mdash dimensional
manifolds to the category of fibred manifolds.
\begin{definition}
A \emph{natural bundle of order~$k$} over~$M$ is any such functor\textup.
\end{definition} 

We shall now give some important examples of (gauge\nobreakdash-) natural bundles.
\begin{example}[Bundle of tensor densities]\label{ex:btd}
A first fundamental example of a natural bundle is given, of course, by the bundle $\Tw^r_sM$ of tensor
densities of weight~$w$ over an \mdash dimensional manifold~$M$. Indeed, $\Tw^r_sM$ is a vector bundle
associated with $L^1\!M$ via the following left action of $G^1_m\cong W^{1,0}_m\{e\}$ on the vector space 
$T^r_s(\R^m)$:
\begin{equation*}
\left\{
\begin{aligned}
 {}&\l\colon G^1_m\times T^r_s(\R^m)\to T^r_s(\R^m)  \\
 {}&\l\colon(\alpha^j\ida_k, t^{p_1\dots p_r}_{q_1\dots q_s})\mapsto
    \alpha^{p_1}\ida_{k_1}\dotsm\alpha^{p_r}\ida_{k_r}t^{k_1\dots k_r}_{l_1\dots l_s}
    \tilde\alpha_{q_1}\ida^{l_1}\dotsm\tilde\alpha_{q_s}\ida^{l_s}(\det\alpha)^{-w}
\end{aligned}
\right.,
\end{equation*}
the tilde over a symbol denoting matrix inversion. For $w=0$ we recover the bundle of tensor fields
over~$M$. This is a definition of $\Tw^r_sM$ which is appropriate for physical applications, where one
usually considers \emph{only} those (active) transformations of tensor fields that are \emph{naturally}
induced by some transformations on the base manifold. Somewhat more unconventionally, though, we can regard
$\Tw^r_sM$ as a \emph{gauge-}natural vector bundle associated with $W^{0,0}(LM)$. Of course, the two bundles
under consideration are the same \emph{as objects}, but their \emph{morphisms} are different.
\end{example}
\begin{example}[Bundle of \Gdash invariant vector fields]
Let $\mathcal{V}:=\R^m \oplus \g$, $\g$ denoting the Lie algebra of~$G$, and consider the following action:
\begin{equation}
\left\{
\begin{aligned}
 {}&\l\colon W^{1,1}_mG\times \mathcal{V}\to \mathcal{V}  \\
 {}&\l\colon\big((\alpha^j\ida_k,a^q,a^r\ida_l),(\nu^i, v^p) \big)\mapsto
 \big(\alpha^i\ida_j\nu^j, A^p\ida_q(a)(v^q +a^q\ida_j\nu^j)\big)
\end{aligned}
\right.,
\label{eq:tau-iv}
\end{equation}
where $(a^q,a^r\ida_l)$ denote natural coordinates on $T^1_mG$: a generic element $j^1_0f\in T^1_mG$ is
represented by $a= f (0)\in G$, i.e.\ $a^q= f^q(0)$, and $a^r\ida_l=\bigl(\de_l(a^{-1}\cdot
f(x))\rvert_{x=0}\bigr)^r$. Obviously, $W^{1,1}\!P\times_\l\mathcal{V}\cong T\!P/G$, its sections thus
representing \Gdash invariant vector fields on~$P$.
\end{example}
\begin{example}[Bundle of vertical \Gdash invariant vector fields]\label{ex:Givfs}
Take~$\g$ as the standard fibre 
and consider the following action:
\begin{equation}
\left\{
\begin{aligned}
 {}&\l\colon W^{1,1}_mG\times\g\to\g  \\
 {}&\l\colon\big((\alpha^j\ida_k,a^q,a^r\ida_l), v^p \big)\mapsto
    A^p\ida_q(a)v^q
\end{aligned}
\right..
\label{eq:tau-viv}
\end{equation}
It is easy to realize that
$W^{1,1}\!P\times_\l\g
\cong V\!P/G\cong (P\times\g)/G$, the bundle of 
vertical \Gdash invariant vector fields on~$P$. Of course, in this example,
we see that~$\g$ is already a $G$-manifold and so $(P\times\g)/G$
is a gauge-natural bundle of order $(0,0)$, i.e.\ a (vector) bundle associated with
$W^{0,0}\!P\cong P$. In other words, giving action~\eqref{eq:tau-viv} amounts to regarding
the original $G$-manifold~$\g$ as a $W^{1,1}_mG$-manifold
via the canonical projection of Lie groups
$W^{1,1}_mG\to G$. It is also meaningful 
to think of action~\eqref{eq:tau-viv} as setting 
$\nu^i=0$ in~\eqref{eq:tau-iv}, and hence one sees 
that the first jet contribution, i.e.~$a^p\ida_i$, disappears.
\end{example}

\section{Split structures on principal bundles}\label{sec:ss}
It is known that, given a principal bundle $P(M,G)$, a \emph{principal connection} on~$P$ may be viewed as
a fibre $G$-equivariant projection $\Phi\colon T\!P\to V\!P$, i.e.\ as a $1$-form in $\Omega^{1}(P,T\!P)$
such that $\Phi\circ\Phi=\Phi$ and $\im\Phi=V\!P$. Here, ``$G$-equivariant'' means that $Tr^a\circ\Phi
=\Phi\circ Tr^a$ for all $a\in G$.  Then, $H\!P:=\ker\Phi$ is a constant-rank vector subbundle of $T\!P$,
called the \emph{horizontal bundle}. We have a decomposition $T\!P=H\!P\oplus V\!P$ and $T_uP=H_uP\oplus
V_uP$ for all $u\in P$. The projection~$\Phi$ is called the \emph{vertical projection} and the  projection
$\chi :=\id_{T\!P} - \Phi$, which is also \Gdash equivariant and satisfies $\chi\circ\chi=\chi$ and $\im\chi
=\ker\Phi$, is called the \emph{horizontal projection}. 

This is, of course, a well-known example of a ``split structure'' on a principal bundle.
We shall now give the following general definition, due---for pseudo-Riemannian manifolds---to a
number of authors \cite{walker55,walker58,cg63,gray67,fava68} and more generally to Gladush and
Konoplya \cite{gk99}.
\begin{definition}
An \emph{$r$-split structure} on a principal bundle $P(M,G)$
is a system of~$r$ fibre \Gdash equivariant linear operators 
$\{\Phi^i\in\Omega^{1}(P,T\!P)\}$\textup, $i=1,2,\dots, r$\textup,
of constant rank with the properties\textup:
\begin{equation}
\Phi^i\cdot\Phi^j =\delta^{ij}\Phi^{j}, \qquad \sum_{i=1}^r\Phi^i =\id_{T\!P}.
\label{eq:Phia}
\end{equation}
\end{definition}
We introduce the notations:
\begin{equation}
\Sigma_u^i :=\im\Phi^i_u,\qquad
n_i :=\dim\Sigma_u^i,
\label{eq:Sigmaa}
\end{equation}
where $\im\Phi^i_u$ is the image of the operator~$\Phi^i$ at a point~$u$ of~$P$,
i.e.\ $\Sigma_u^i=\{\,v\in T_uP\mid\Phi^i_u\circ v=v\,\}$.
Owing to the constancy of the rank of the operators~$\{\Phi^i\}$,
the numbers~$\{n_i\}$ do not depend on the point~$u$ of~$P$.
It follows from the very definition of an $r$-split structure
that we have a \Gdash equivariant decomposition of the
tangent space:
\begin{equation*}
T_uP = {\bigoplus_{i=1}^{r}}\Sigma_u^i,\qquad
\dim T_uP =\sum_{i=1}^r n_i.
\end{equation*}
Obviously, the bundle $T\!P$ is also decomposed into $r$ vector
subbundles~$\{\Sigma^i\}$ so that
\begin{equation}
T\!P ={\bigoplus_{i=1}^{r}}\Sigma^i,\qquad
\Sigma^i ={\bigcup_{u\in P}}\Sigma_u^i.
\label{eq:TPdecomp}
\end{equation}
\begin{remark}
In general, the $r$ vector subbundles $\{\Sigma^i\to P\}$ 
are \emph{anholonomic}, i.e.\ non-integrable, and are not vector subbundles of $V\!P$.
For a principal connection, i.e.\ for the case $T\!P=H\!P\oplus V\!P$,
the subbundle $V\!P$ is integrable.
\end{remark}
\begin{proposition}
An equivariant decomposition of $T\!P$ into $r$ vector
subbundles $\{\Sigma^i\}$ as given by~\eqref{eq:TPdecomp}\textup, with 
$T_ur^a(\Sigma_u^i)=\Sigma_{u\cdot a}^i$\textup, induces 
a system of~$r$ fibre \Gdash equivariant linear operators 
$\{\Phi^i\in \Omega^{1}(P,T\!P)\}$ of constant rank 
satisfying properties~\eqref{eq:Phia} and~\eqref{eq:Sigmaa}\textup.
\end{proposition}
\begin{proposition}
Given an \emph{$r$-split structure} on a principal bundle $P(M,G)$\textup,
every \Gdash invariant vector field~$\Xi$ on~$P$ splits into
$r$ invariant vector fields~$\{\Xi_i\}$ such that 
$\Xi=\Xi_1\oplus\dots\oplus\Xi_r$ and $\Xi_i(u)\in\Sigma^i_u$ 
for all $u\in P$ and $i=1,2,\dots, r$\textup.
\end{proposition}
\begin{remark}
The vector fields $\{\Xi_i\}$ are compatible with the $\{\Sigma^i\}$, i.e.\ they are sections 
$\{\Xi_i\colon P\to\Sigma^i\}$ of the vector bundles $\{\Sigma^i\to P\}$.
\end{remark}
\begin{corollary}
Let $P(M,G)$ be a reductive \Gdash structure
on a principal bundle $Q(M,H)$ and let $i_P\colon P\to Q$ be the canonical
embedding\textup. Then\textup, any given $r$-split structure 
on $Q(M,H)$ induces an $r$-split structure restricted to $P(M,G)$\textup, 
i\textup.e\textup.\ an equivariant decomposition of 
$i_{P}^{\ast}(TQ)\equiv P\times_{Q}TQ=
\{\,(u,v)\in P\times TQ \mid i_{P}(u)=\tau_{Q}(v)\,\}$
such that $i_{P}^{\ast}(TQ)=i_{P}^{\ast}(\Sigma^{1})
\oplus\dots\oplus i_{P}^{\ast}(\Sigma^{r})$\textup,
and any \Hdash invariant vector field~$\Xi$ on~$Q$
restricted to~$P$ splits into $r$ \Gdash invariant sections
of the pull-back bundles $\{i_{P}^{\ast}(\Sigma^i)\equiv P\times_{Q}\Sigma^i\}$\textup,
i\textup.e\textup.\ $\Xi=\Xi_1\oplus\dots\oplus\Xi_r$ with $\Xi_i(u)\in\bigl(i_P^\ast(\Sigma^i)\bigr)_u$ 
for all $u\in P$ and $i\in\{1,2,\dots, r\}$\textup.
\end{corollary}
\begin{remark}\label{rem:Hinv}
Note that the pull-back $i_{P}^{\ast}$ is a \emph{natural operation}, i.e.\ it respects the splitting
$i_{P}^{\ast}(TQ)=i_{P}^{\ast}(\Sigma^{1})\oplus\dots\oplus i_{P}^{\ast}(\Sigma^{r})$. In other words, the
pull-back of a splitting for~$Q$ is a splitting of the pull-backs for~$P$. Furthermore, although the vector
fields~$\{\Xi_i\}$ are \Gdash invariant sections of their respective pull-back bundles, they are \Hdash
invariant if regarded as vector fields on the corresponding subsets of~$Q$.
\end{remark}
In \S\ref{sec:gnb} we saw that $W^{k,h}\!P$ is a principal bundle over $M$.
Consider in particular $W^{1,1}\!P$, the $(1,1)$-principal prolongation of $P$. 
The fibred manifold $W^{1,1}\!P\to M$ coincides with the fibred product 
$W^{1,1}\!P:= L^1\!M\times_{M} J^1\!P$ over~$M$.
We have two canonical principal bundle morphisms $\pr_1\colon W^{1,1}\!P\to L^1\!M$
and $\pr_2\colon W^{1,1}\!P\to P$. In particular, $\pr_2\colon W^{1,1}\!P\to P$ is a
$G^1_m\rtimes\g\otimes \R^{m}$-principal bundle, $G^1_m \rtimes \g\otimes \R^{m}$ being the
kernel of $W^{1}_{m}G\to G$\textup. The following lemma recognizes $\tau_{P}\colon T\!P \to P$
as a vector bundle associated with the principal bundle $W^{1,1}\!P\to P$\textup.
\begin{lemma}
The vector bundle $\tau_{P}\colon T\!P \to P$ coincides with
the vector bundle $T^{1,1}\!P := (W^{1,1}\!P\times \mathcal{V})/(G^1_m\rtimes\g\otimes\R^m)$
over $P$\textup, where $\mathcal{V}:=\R^m\oplus \g$ is the left $G^1_m\rtimes\g\otimes\R^m$-manifold
with action given by\textup:
\begin{equation}
\left\{
\begin{aligned}
 {}&\tau\colon G^1_m \rtimes \g\otimes \R^{m}\times \mathcal{V}\to \mathcal{V}  \\
 {}&\tau\colon\big((\alpha^j\ida_k,e,a^r\ida_l),(\nu^i, v^p) \big)\mapsto
 (\alpha^i\ida_j\nu^j, v^p + a^p\ida_i\nu^i)
\end{aligned}
\right..
\label{eq:tau-vv}
\end{equation}
\end{lemma}
\begin{remark}
The vector bundle $\tau_{P}\colon T\!P \to P$
is a gauge-natural bundle of order $(0,0)$ associated with the 
$G^1_m \rtimes \g\otimes \R^{m}$-principal bundle
$\pr_2\colon W^{1,1}\!P\to P$.
\end{remark}
\begin{lemma}
$V\!P \to P$ is a trivial vector bundle associated with $W^{1,1}\!P\to P$\textup.
\end{lemma}
\begin{lemma}
Let $P(M,G)$ be a reductive \Gdash structure
on a principal bundle $Q(M,H)$ and $i_P\colon P\to Q$ the canonical
embedding\textup. Then\textup, $i_P^{\ast}(TQ) = P\times_{Q}TQ$ is a vector bundle 
over~$P$ associated with $W^{1,1}P\!\to P$\textup.
\end{lemma}
From the above lemmas it follows that
another important example of a split structure on a principal bundle
is given by the following
\begin{theorem}\label{thm:ft}
Let $P(M,G)$ be a reductive \Gdash structure
on a principal bundle $Q(M,H)$ and let $i_P\colon P\to Q$ be the canonical
embedding\textup. Then\textup, there exists 
a canonical decomposition of $\, i_P^{\ast}(TQ)\to P$ such that
\begin{equation*}
i_{P}^{\ast}(TQ) = T\!P \oplus\mathcal M(P),
\end{equation*}
i\textup.e\textup.\ at each $u\in P$ one has
\begin{equation*}
T_uQ = T_uP\oplus\mathcal M_u,
\end{equation*}
$\mathcal M_u$ being the fibre over~$u$ of the subbundle $\mathcal M(P)\to P$
of $i_P^{\ast}(V\!Q)\to P$\textup. The bundle $\mathcal M(P)$ is defined as
$\mathcal M(P) := (W^{1,1}\!P\times \m)/(G^1_m \rtimes\g\otimes\R^{m})$\textup, where $\m$ 
is the \textup(trivial left\textup) $G^1_m\rtimes\g\otimes\R^{m}$-manifold\textup. 
\end{theorem}
\begin{remark}
The trivial $G^1_m \rtimes \g\otimes \R^{m}$-manifold~$\m$ corresponds to the
action~\eqref{eq:tau-viv} of Example~\ref{ex:Givfs} with $W^{1,1}_{m}G$
restricted to~$G^1_m \rtimes \g\otimes \R^{m}$, and $\g$ restricted to~$\m$.
Of course, since the group $G^1_m \rtimes \g\otimes \R^{m}$ acts trivially on~$\m$,
it follows that $\mathcal M(P)$ is trivial, i.e.\ isomorphic to
$P\times \m$, because $W^{1,1}\!P/(G^1_m \rtimes \g\otimes \R^{m})\cong P$.
\end{remark}
From the above theorem two corollaries follow, which are of prime importance for the concepts of a Lie
derivative we shall introduce in the next section.
\begin{corollary}\label{cor:gKosmann}
Let $P(M,G)$ and $Q(M,H)$ be as in the previous theorem\textup. The restriction~$\Xi\rvert_P$ of an \Hdash 
invariant vector field~$\Xi$ on~$Q$ to~$P$ splits into a \Gdash invariant vector field~$\Xi_{\mathrm K}$
on~$P$\textup, called the \emph{generalized Kosmann vector field associated with~$\Xi$}\textup, and a
``transverse'' vector field~$\Xi_{\mathrm G}$, called the \emph{generalized von G\"oden vector field
associated with~$\Xi$}\textup.
\end{corollary}
\begin{corollary}\label{cor:Kosmann}
Let $P(M,G)$ be a classical \Gdash structure\textup, i\textup.e\textup.\ a reductive \Gdash structure
on the bundle~$LM$ of linear frames over~$M$\textup. The restriction~$L\xi\rvert_P$ to $P \to M$ of 
the natural lift $L\xi$ onto~$ LM$ of a vector field~$\xi$ on~$M$  splits into a \Gdash invariant vector
field on~$P$ called the \emph{generalized Kosmann lift of~$\xi$} and denoted simply by $\xi_{\mathrm
K}$\textup, and a ``transverse'' vector field called the \emph{von G\"oden lift of~$\xi$} and denoted by
$\xi_{\mathrm G}$\textup.
\end{corollary}
\begin{remark}
The last corollary still holds if, instead of $ LM$, one considers the $k$-th order frame bundle $L^k\!M$ 
and hence a classical \Gdash structure of order $k$, i.e.\ a reductive \Gdash subbundle $P$ of $L^k\!M$.
Note also that the Kosmann lift $\xi\mapsto\xi_{\mathrm K}$ is \emph{not} a Lie algebra homomorphism,
although $\xi_{\mathrm K}$ is a \Gdash invariant vector field and projects over~$\xi$.
\end{remark}
\begin{example}[Kosmann lift]\label{ex:Kl}
A fundamental example of a \Gdash structure on a manifold~$M$ is given, of course, by the
bundle~$\SO(M,g)$ of its (pseudo\nobreakdash-) orthonormal frames with respect to a metric~$g$ of signature
$(p,q)$, where $p+q=m\equiv\dim M$. $\SO(M,g)$ is a principal bundle (over~$M$) with structure
group $G=\SO(p,q)$. Now, recall that the natural lift of a vector field~$\xi$ onto~$LM$ is defined as
\begin{equation*}
L\xi:=\left.\frac\de{\de t}L^1\varphi_t\right\rvert_{t=0},
\end{equation*}
$\{\varphi_t\}$ denoting the flow of~$\xi$. If $(\rho_a\ida^b)$ denotes a
(local) basis of right \GLdash invariant vector fields on $LM$ reading $(\rho_a\ida^b=u_c\ida^b\de/\de
u_c\ida^a)$ in some local chart $(x^\mu,u_a\ida^b)$ and $(e_a =:e_a\ida^\mu\de_\mu)$ is a local section
of~$LM$, then $L\xi$ has the local expression
\begin{equation*}
L\xi =\xi^a e_a +(L\xi)^a\ida_b\rho_a\ida^b,
\end{equation*}
where $\xi=:\xi^a e_a$ and
\begin{equation*}
(L\xi)^a\ida_b :=\tilde e^a\ida_\rho(\de_\nu\xi^\rho e_b\ida^\nu -\xi^\nu\de_\nu e_b\ida^\rho).
\end{equation*}
If we now let $(e_a)$ and $(x^\mu,u_a\ida^b)$ denote a local section and a local chart of $\SO(M,g)$,
respectively, then  the generalized Kosmann lift~$\xiK$ on $\SO(M,g)$ of a vector field~$\xi$ on~$M$, simply
called its \emph{Kosmann lift} \cite{fffg96}, locally reads
\begin{equation*}
\xiK = \xi^a e_a +(L\xi)_{[ab]}A^{ab},
\end{equation*}
where $(A^{ab})$ is a basis of right \SOdash invariant vector fields on $\SO(M,g)$ locally reading
$(A^{ab}=\eta^{c[a}\delta^{b]}\ida_d\rho_c\ida^d)$, $(L\xi)_{ab} :=\eta_{ac}(L\xi)^c\ida_b$, and
$(\eta_{ac})$ denote the components of the standard Minkowski metric of  signature $(p,q)$.
\end{example}
Now, combining Proposition~\ref{prop:GivfGaivf} and Theorem~\ref{thm:ft} yields the following result, which,
in particular, will enable us to extend the concept of a Kosmann lift to the important context of spinor
fields.
\begin{corollary}\label{cor:GaxiK}
Let $\zeta \colon \tilde P \to P$ be a \Gadash structure 
over a classical \Gdash structure $P(M,G)$\textup. 
Then\textup, the generalized Kosmann lift
$\xi_{\mathrm K}$ of a vector field~$\xi$ on~$M$ lifts 
to a unique \textup(\Gadash invariant\textup) vector
field~$\tilde\xi_{\mathrm K}$ on~$\tilde P$\textup, 
which projects over~$\xi_{\mathrm K}$\textup.
\end{corollary}

\section{Lie derivatives on reductive \protect\Gdash structures}\label{sec:ld}
As already mentioned in the Introduction, the general theory of Lie derivatives 
stems from Trautman's seminal paper \cite{trautman72}. 
Here, we mainly follow the notation and conventions
of \cite[\S47]{kms}.
\begin{definition}
Let $M$ and $N$ be two manifolds and $f\colon M\to N$ 
a map between them\textup. By a \emph{vector
field along~$f$} we shall mean a map $Z\colon M\to T\!N$ 
such that $\tau_N\circ Z=f$, $\tau_N\colon T\!N\to
N$ denoting the canonical tangent bundle projection.
\end{definition}
\begin{definition}\label{dfn:gLieXY}
Let $M$\textup, $N$ and $f$ be as above\textup, 
and let $X$ and~$Y$ be two vector fields on~$M$
and~$N$\textup, respectively\textup. 
Then\textup, by the \emph{generalized Lie derivative
$\Lie_{(X,Y)}f$ of~$f$ with respect to~$X$ and~$ Y$} 
we shall mean the vector field along~$f$ given by
\begin{equation*}
\Lie_{(X,Y)}f :=Tf\circ X- Y\circ f.
\end{equation*}
\end{definition}
If $\{\varphi_t\}$ and $\{\Phi_t\}$ denote the flows 
of~$X$ and~$Y$, respectively, then one readily verifies
that
\begin{equation*}
\Lie_{(X,Y)}f=\left.\frac\de{\de t}
  (\Phi_{-t}\circ f\circ\varphi_{t})\right\rvert_{t = 0}.
\end{equation*}
%
An important specialization of Definition~\ref{dfn:gLieXY} is given by the following
\begin{definition}\label{dfn:gLieXI}
Let $\pi\colon B\to M$ be a fibred manifold\textup, 
$\sigma\colon M\to B$ a section of~$\pi$\textup, and~$\Xi$
a projectable vector field on~$B$ over 
a vector field~$\xi$ on~$M$\textup. Then\textup, by the
\emph{generalized Lie derivative $\Lie_\Xi\sigma$ 
of~$\sigma$ with respect to~$\Xi$} we shall
mean the map
\begin{equation}
\Lie_\Xi\sigma :=\Lie_{(\xi,\Xi)}\sigma\colon M\to V\!B.
\end{equation}
\end{definition}
(It is easy to realize that 
$\Lie_{\Xi}\sigma\equiv T\sigma\circ\xi -\Xi\circ\sigma$ 
takes indeed values in the vertical tangent bundle simply by applying 
$T\pi$ to it and remembering that $\Xi$ is projectable.)

Now recall that a fibred manifold $\pi\colon B \to M$
admits a \emph{vertical splitting} if there exists a linear bundle isomorphism 
(covering the identity of~$B$) 
$\alpha \colon V\!B \to B\times_{M}{\bar B}$, 
where ${\bar\pi}\colon {\bar B} \to M$ is a 
vector bundle. In particular, a vector bundle $\pi\colon B \to M$ 
admits a \emph{canonical} vertical
splitting $\alpha\colon V\!B \to B\times_{M}B$. Indeed, if  $\tau_B\colon T\!B \to B$ 
denotes the (canonical) tangent bundle projection restricted to $V\!B$, $y$ is a point in~$B$ such that 
$y=\tau_B(v)$ for a given $v\in V\!B$, and $\gamma\colon\R \to B_y\equiv\pi^{-1}\bigl(\pi(y)\bigr)$ is a
curve such that $\gamma(0)=y$ and $j^1_0\gamma=v$, then~$\alpha$ is given by $\alpha (v) := (y, w)$, where
$w :=\lim_{t \to 0}\frac1t{( \gamma(t) - \gamma(0) )}$.

\begin{proposition}\label{prop:spLie} 
In this case\textup, the generalized Lie derivative
$\Lie_\Xi\sigma$ is of the form 
\begin{equation}
\Lie_\Xi\sigma = (\sigma,\,£_\Xi\sigma ),
\end{equation}
the first component being the original section $\sigma$\textup.
The second component $£_{\Xi}\sigma$ is a section of~$\bar B$\textup,
called the \emph{Lie derivative of $\sigma$ with respect to $\Xi$}\textup.
For the sake of clarity\textup, the operator~$\pounds$ will be occasionally referred to as the
\emph{restricted Lie derivative} \textup{\cite[\S47]{kms}.}
\end{proposition}
\begin{remark}
In this case, on using the fact that $£_\Xi\sigma$ is the derivative of
$\Phi_{-t}\circ\sigma\circ\varphi_{t}$ at $t =0$ in the classical sense, 
one can re-express the restricted Lie derivative in the form
\begin{equation}
£_\Xi\sigma(x) =\lim_{t \to 0}\frac1t\bigl(\Phi_{-t}\circ\sigma\circ\varphi_{t}(x) -\sigma(x)\bigr).
\end{equation}
\end{remark}
Proposition~\ref{prop:spLie} also works whenever $B$ is an \emph{affine} bundle.
This is so because, also in this case, $\pi\colon B \to M$ admits a canonical vertical
decomposition $\alpha \colon V\!B \to B\times_{M}{\bar B}$, 
where ${\bar\pi}\colon {\bar B} \to M$ is the
vector bundle  associated with~$B$.

Now, we can specialize Definition~\ref{dfn:gLieXI} 
to the case of gauge-natural bundles in a straightforward
manner.
\begin{definition}\label{def:gnld}
Let $P_\l$ be a gauge-natural bundle associated with some principal bundle $P(M,G)$, $\Xi$ a \Gdash 
invariant vector field on~$P$ projecting over a vector field~$\xi$ on~$M$\textup, and $\sigma\colon M\to
P_\l$ a section of~$P_\l$\textup. Then, by the \emph{generalized \textup(gauge-natural\textup) Lie derivative
of~$\sigma$ with respect to~$\Xi$} we shall mean the map
\begin{equation}
\Lie_\Xi\sigma\colon M\to V\!P_\l, \quad
\Lie_\Xi\sigma :=T\sigma\circ\xi -\Xi_\l\circ\sigma,
\label{eq:gnld}
\end{equation}
where $\Xi_\l$ is the generator of the \onedash 
parameter group~$\{(\Phi_t)_\l\}$ of automorphisms of~$P_\l$
functorially induced by the flow~$\{\Phi_t\}$ of~$\Xi$ \textup[\cf~\eqref{eq:indaut}\textup{].}
Equivalently\textup,
\begin{equation}
\Lie_\Xi\sigma =\left.\frac\de{\de t}
  \bigl((\Phi_{-t})_\l\circ \sigma\circ\varphi_t\bigr)\right\rvert_{t = 0},
\tag{$\ref{eq:gnld}'$}
\end{equation}
$\{\varphi_t\}$ denoting the flow of~$\xi$\textup.
\end{definition}
As usual, whenever $\pi\colon P_\l\to M$ admits 
a canonical vertical splitting of $V\!P_\l$,
we shall write $£_\Xi\sigma\colon M\to\bar P_\l :=\overline{P_\l}$ 
for the corresponding restricted Lie derivative.

Furthermore, for each \Gadash structure 
$\zeta \colon \tilde P \to P$ on $P$,  we shall simply write
$£_\Xi\tilde\sigma :=£_{\tilde\Xi}\tilde\sigma\colon M\to\bar{\tilde P}_{\tilde\l}$, $\tilde P_{\tilde\l}$
denoting a gauge-natural bundle associated with~$\tilde P$ 
(admitting a canonical vertical splitting) and
$\tilde\sigma\colon M\to\tilde P_{\tilde\l}$ one of its sections,  
since~$\Xi$ admits a unique (\Gadash invariant) 
lift $\tilde\Xi$ onto~$\tilde P$ (\cf\ Proposition~\ref{prop:GivfGaivf}).
%
%
We stress that \emph{Definition~\ref{def:gnld} is the conceptually natural generalization of the
classical notion of a Lie derivative} \cite{yano57}, to which it suitably reduces when applied to natural
objects and, hence, notably, to tensor fields and tensor densities.

Of course, we can now further specialize to the case of classical \Gdash structures and, in particular,
give the following
\begin{definition}\label{def:ldK}
Let $P_\l$ be a gauge-natural bundle associated with some classical \Gdash structure $P(M,G)$\textup,
$\xi_{\mathrm K}$ the generalized Kosmann lift \textup(on~$P$\textup) of a vector field~$\xi$
on~$M$\textup,  and $\sigma\colon M\to P_\l$ a section of~$P_\l$\textup. Then\textup, by the
\emph{generalized Lie derivative $\Lie_\xi\sigma$  of~$\sigma$ with respect to~$\xi$} we shall mean the map
$\Lie_\xi\sigma :=\Lie_{\xi_{\mathrm K}}\sigma$\textup,  where $\Lie_{\xi_{\mathrm K}}\sigma$ denotes 
the generalized Lie derivative of~$\sigma$ with respect to~$\xi_{\mathrm K}$ in the sense of
Definition~\ref{def:gnld}\textup.
\end{definition}
Consistently, we shall simply write 
$£_\xi\sigma :=£_{\xi_{\mathrm K}}\sigma\colon M\to\bar P_\l$ for the 
corresponding restricted Lie derivative, whenever defined, 
and $£_\xi\tilde\sigma :=£_{\tilde\xi_{\mathrm K}}\tilde\sigma\colon M\to\bar{\tilde
P}_{\tilde\l}$ for the (restricted) 
Lie derivative of a section~$\sigma$ of a gauge-natural bundle~$\tilde
P_{\tilde\l}$ associated with some principal 
prolongation of a \Gadash structure $\zeta\colon\tilde P\to P$
(and admitting a canonical vertical splitting), 
which makes sense since $\xi_{\mathrm K}$ admits a unique
(\Gadash invariant) lift $\tilde\xi_{\mathrm K}$ 
onto~$\tilde P$ (\cf\ Corollary~\ref{cor:GaxiK}).

\begin{example}[Lie derivative of spinor fields. I]
In Example~\ref{ex:Kl} we mentioned that a fundamental example of a \Gdash structure on a
manifold~$M$ is given by the bundle~$\SO(M,g)$ of its (pseudo\nobreakdash-) orthonormal frames. An equally
fundamental example of a \Gadash structure on $\SO(M,g)$ is given by the corresponding spin
bundle $\Spin(M,g)$ with structure group $\Gamma=\Spin(p,q)$. 
Now, it is obvious that spinor fields can be 
regarded as sections of a suitable gauge-natural
bundle over~$M$. Indeed, if $\l$ is 
the linear representation of $\Spin(p,q)$ on the vector
space $\C^m$ induced by a given choice of $\gamma$ matrices, 
then the associated vector bundle
$S(M) :=\Spin(M,g)\times_\l\C^m$ is 
a gauge-natural bundle of order $(0,0)$ whose sections
represent spinor fields (or, more precisely, 
spin-vector fields). Therefore, in spite of
what is sometimes believed, 
a Lie derivative of spinors (in the sense of
Definition~\ref{def:gnld}) always exists, 
\emph{no matter what} the vector field~$\xi$ on~$M$ is. 
Locally, such a Lie derivative reads
\begin{equation*}
£_\Xi\psi =\xi^a e_a\psi +\frac14\Xi_{ab}\ga^a\ga^b\psi
\end{equation*}
for any spinor field~$\psi$, 
$(\Xi_{ab}=\Xi_{[ab]})$ denoting
the components of an \SOdash invariant 
vector field $\Xi=\xi^a e_a +\Xi_{ab}A^{ab}$
on $\SO(M,g)$, $\xi=:\xi^a e_a$, and $e_a\psi$ the Pfaff derivative of~$\psi$ along the local section
$(e_a=:e_a\ida^\mu\de_\mu)$ of $\SO(M,g)$ induced by some local section of $\Spin(M,g)$. This is the most
general  notion of a (gauge-natural) Lie derivative of spinor fields and the appropriate one for most
situations of physical interest (\cf\ \cite{gmv,gng}): the generality of~$\Xi$ might be disturbing, but is
the \emph{unavoidable} indication that $S(M)$ is \emph{not} a natural bundle.

\noindent
If we wish nonetheless to remove such a generality, 
we must \emph{choose} some canonical
(\emph{not} natural) lift of~$\xi$ onto~$\SO(M,g)$. 
The conceptually (\emph{not} mathematically)
most ``natural'' choice is perhaps given 
by the Kosmann lift (recall Example~6 and use
Corollary~\ref{cor:GaxiK}). The ensuing Lie derivative locally reads
\begin{equation}
£_\xi\psi =\xi^a e_a\psi +\frac14(L\xi)_{[ab]}\ga^a\ga^b\psi.
\label{eq:dlK}
\end{equation}
Of course, if `$\nabla$' denotes the covariant derivative operator associated with the Levi-Civita (or
Riemannian) connection with respect to~$g$, the previous expression can be recast into the form
\begin{equation}
£_\xi\psi =\xi^a\cde_a\psi -\frac14\cde_{[a}\xi_{b]}\ga^a\ga^b\psi,
\tag{$\ref{eq:dlK}'$}
\label{eq:dlKp}
\end{equation}
which reproduces exactly Kosmann's definition \cite{kosmann72} (see \cite{fffg96} for further details and a
more thorough discussion).  We stress that, although \emph{in this case} its local expression would be
identical  with~\eqref{eq:dlK}, this is \emph{not} the ``metric Lie derivative'' introduced by Bourguignon
and Gauduchon in~\cite{bg92}. To convince oneself of this it is enough to take the Lie derivative of the
metric~$g$, which is a section of the \emph{natural} bundle~$\Symm^2T^*\!M$, `$\Symm$' denoting the
symmetrized tensor product. Since the (restricted) Lie derivative~$£_\xi$ in the sense of
Definition~\ref{def:ldK} must reduce to the ordinary one on natural objects, it holds that 
\begin{equation*}
£_{L\xi}g =£_\xi g.
\end{equation*}
On the other hand, if $£_\xi$ coincided 
with the operator~$£_\xi^g$ defined by Bourguignon and
Gauduchon,  the right-hand side of the above identity should equal zero
\cite[Proposition~15]{bg92}, thereby 
implying that $\xi$ is a Killing vector field, contrary to
the fact that~$\xi$ is completely arbitrary. 
Indeed, in order to recover Bourguignon and
Gauduchon's definition, another concept 
of a Lie derivative must be introduced.
\end{example}
%
%
We shall start by recalling two
classical definitions \cite{koba72}.
\begin{definition}
Let $P(M,G)$ be a \textup(classical\textup) \Gdash structure\textup. Let $\varphi$ be a
diffeomorphism of~$M$ onto itself and $L^1\varphi$ its natural lift onto~$LM$\textup. If
$L^1\varphi$ maps $P$ onto itself\textup, i\textup.e\textup.\ if $L^1\varphi(P)\subseteq P$\textup, then
$\varphi$ is called an \emph{automorphism} of the \Gdash structure~$P$\textup.
\end{definition}
\begin{definition}
Let $P(M,G)$ be a \Gdash structure\textup. A vector field~$\xi$ on~$M$ is called an
\emph{infinitesimal automorphism} of the \Gdash structure~$P$ if it generates a local
\onedash parameter group of automorphisms of~$P$\textup.
\end{definition}
We can now generalize these concepts to the framework of reductive \Gdash structures as follows.
\begin{definition}
Let $P(M,G)$ be a reductive \Gdash structure on a principal bundle $Q(M,H)$ and $\Phi$ a principal
automorphism of~$Q$\textup. If $\Phi$ maps $P$ onto itself\textup, i\textup.e\textup.\ if $\Phi(P)\subseteq
P$\textup, then $\Phi$ is called a \emph{generalized automorphism} of the reductive \Gdash
structure~$P$\textup.
\end{definition}
Of course, each element of $\Aut(P)$, i.e.\ each principal automorphism of~$P$, is by definition
a generalized automorphism of the reductive \Gdash structure~$P$. Analogously, we have
\begin{definition}
Let $P(M,G)$ be a reductive \Gdash structure on a principal bundle $Q(M,H)$\textup. An \Hdash invariant
vector field~$\Xi$ on~$Q$ is called a \emph{generalized infinitesimal automorphism} of the reductive \Gdash
structure~$P$ if it generates a local \onedash parameter group of generalized automorphisms of~$P$\textup.
\end{definition}
Of course, each element of $\X_G(P)$, i.e.\ each \Gdash invariant vector field on~$P$, is by
definition a generalized infinitesimal automorphism of the reductive \Gdash structure~$P$.

Now, along the lines of \cite[Proposition X.1.1]{kn2} it is easy to prove
\begin{proposition}
Let $P(M,G)$ be a reductive \Gdash structure on a principal bundle $Q(M,H)$\textup. An \Hdash invariant
vector field~$\Xi$ on~$Q$ is a generalized infinitesimal automorphism of the reductive \Gdash structure~$P$
if and only if~$\Xi$ is tangent to $P$ at each point of $P$\textup.
\end{proposition}
We then have the following important
\begin{lemma}
Let $P(M,G)$ be a reductive \Gdash structure on a principal bundle $Q(M,H)$ and~$\Xi$ a generalized
infinitesimal automorphism of the reductive \Gdash structure~$P$\textup. Then\textup, the flow $\{\Phi_t\}$
of~$\Xi$, it being \Hdash invariant, induces on each gauge-natural bundle $Q_\l$ associated with~$Q$ a
\onedash parameter group $\{(\Phi_t)_\l\}$ of global automorphisms\textup.
\end{lemma}
\begin{proof}
Since $\Xi$ is by assumption a generalized infinitesimal automorphism, it is by definition an \Hdash
invariant vector field on~$Q$. Therefore, its flow $\{\Phi_t\}$ is a \onedash parameter group of \Hdash
equivariant maps on~$Q$. Then, if $Q_\l=W^{k,h}Q\times_\l F$, we set
\begin{equation*}
(\Phi_t)_\l([u, f]_\l) := [W^{k,h}\Phi_t(u),f]_\l,
\end{equation*}
$u\in Q$, $f\in F$, and are back to the situation of formula~\eqref{eq:indaut}.
\end{proof}
%
%
\begin{corollary}
Let $P(M,G)$ and $Q(M,H)$ be as in the previous lemma\textup, and let~$\Xi$ be an \Hdash invariant
vector field on~$Q$\textup. Then\textup, the flow $\{(\Phi_{\mathrm K})_t\}$ of the generalized Kosmann
vector field~$\Xi_{\mathrm K}$ associated with~$\Xi$ induces on each gauge-natural bundle $Q_\l$ associated
with~$Q$ a \onedash parameter group $\bigl\{\bigl((\Phi_{\mathrm K})_t\bigr)_\l\bigr\}$ of global
automorphisms\textup.
\end{corollary}
\begin{proof}
Recall that, although the generalized Kosmann vector field~$\Xi_{\mathrm K}$ is a \Gdash invariant vector
field on~$P$, it is \Hdash invariant if regarded as a vector field on the corresponding subset of~$Q$
(\cf\ Remark~\ref{rem:Hinv} and Corollary~\ref{cor:gKosmann}). Therefore, its flow $\{(\Phi_{\mathrm
K})_t\}$ is a \onedash parameter group of \Hdash equivariant automorphisms on the subset~$P$ of~$Q$.

\noindent
We now want to define a \onedash parameter group of automorphisms $\bigl\{\bigl((\Phi_{\mathrm
K})_t\bigr)_\l\bigr\}$ of $Q_\l=W^{k,h}Q\times_\l F$. Let $[u, f]_\l\in Q_\l$, $u\in Q$ and $f\in F$, and let
$u_1$ be a point in~$P$ such that $\pi(u_1)=\pi(u)$, $\pi\colon Q\to M$ denoting the canonical projection.
There exists a unique
$a_1\in H$ such that
$u = u_1\cdot a_1$. Set then
\begin{equation*}
\bigl((\Phi_{\mathrm K})_t\bigr)_\l([u, f]_\l) :=[W^{k,h}(\Phi_{\mathrm K})_t(u_1),
a_1\cdot f]_\l.  
\end{equation*}
We must show that, given another point $u_2\in P$ such that $u = u_2\cdot a_2$ for some (unique) $a_2\in H$,
we have
\begin{equation*}
[W^{k,h}(\Phi_{\mathrm K})_t(u_1), a_1f]_\l = [W^{k,h}(\Phi_{\mathrm K})_t(u_2), a_2\cdot f]_\l.
\end{equation*}
Indeed, since the action of~$H$ is free and transitive on the fibres, from $u = u_1\cdot a_1$ and $u =
u_2\cdot a_2$ it follows that $a_1 = a\cdot a_2$ or $a = a_1\cdot (a_2)^{-1} $ or $a_2 = a^{-1}\cdot a_1$.
But then
\begin{align*}
[W^{k,h}(\Phi_{\mathrm K})_t(u_2), a_2\cdot f]_\l
&= [W^{k,h}(\Phi_{\mathrm K})_t(u_1\cdot a), a^{-1}\cdot a_1\cdot f]_\l   \\
&= [W^{k,h}(\Phi_{\mathrm K})_t(u_1)\odot W^{k,h}_ma, a^{-1}\cdot a_1\cdot f]_\l   \\
&= [W^{k,h}(\Phi_{\mathrm K})_t(u_1), a_1\cdot f]_\l,
\end{align*}
as claimed. It is then easy to see that the so-defined $\bigl((\Phi_{\mathrm K})_t\bigr)_\l$ does not depend
on the chosen representative.
\end{proof}
By virtue of the previous corollary, we can now give the following
\begin{definition}\label{def:rld}
Let $P(M,G)$ be a reductive \Gdash structure on a principal bundle $Q(M,H)$\textup, $G\neq\{e\}$\textup, and
$\Xi$ an \Hdash invariant vector field on~$Q$ projecting over a vector field~$\xi$ on~$M$\textup. 
Let $Q_\l$ be a gauge-natural bundle associated with~$Q$ and $\sigma\colon M\to Q_\l$ a section
of~$Q_\l$\textup. Then\textup, by the \emph{generalized \Gdash  reductive Lie derivative of~$\sigma$ with
respect to~$\Xi$}  we shall mean the map
\begin{equation*}
\Lie_\Xi^G\sigma :=\left.\frac\de{\de t}
  \Bigl(\bigl((\Phi_{\mathrm K})_{-t}\bigr)_\l\circ \sigma\circ\varphi_t\Bigr)\right\rvert_{t = 0},
\end{equation*}
$\{\varphi_t\}$ denoting the flow of~$\xi$\textup.
\end{definition}
The corresponding notions of a restricted Lie derivative and a (generalized or restricted)  Lie
derivative on an associated \Gadash structure can be defined in the usual way. 
\begin{remark}
When $Q=P$ (and $H=G$), $\Xi_{\mathrm K}$ is just $\Xi$, and we recover the notion of a (generalized) Lie
derivative in the sense of Definition~\ref{def:gnld}, \emph{but}, as $G$ is required not to equal the trivial
group $\{e\}$, $Q_\l$ is never allowed to be a (purely) natural bundle.
\end{remark}
By its very definition, the (restricted) \Gdash reductive Lie derivative does \emph{not} reduce, in general,
to the ordinary (natural) Lie derivative on fibre bundles associated with $L^k\!M$. This fact makes it
unsuitable in all those situations where one needs a \emph{unique} operator which reproduce ``standard
results'' when applied to ``standard objects''.

In other words, $£_\Xi^G$ is defined with respect to some \emph{pre-assigned} (generalized)
symmetries. We shall make this statement explicit in Proposition~\ref{prop:god00} below, which provides a
generalization of a well-known classical result. 

Let then $\KK$ be a tensor over the vector space $\R^m$ (i.e., an element of the tensor algebra over $\R^m$)
and $G$ the group of linear transformations of $\R^m$ leaving $\KK$ invariant. Recall that each reduction of
the structure group $\GL(m,\R)$ to $G$ gives rise to a tensor field~$K$ on~$M$. Indeed, we may regard each
$u\in LM$ as a linear isomorphism of~$\R^m$ onto $T_xM$, where $x=\pi(u)$ and $\pi\colon LM\to M$ denotes, as
usual, the canonical projection. Now, if $P(M,G)$ is a \Gdash structure, at each  point~$x$ of~$M$ we can
choose a frame~$u$ belonging to $P$ such that $\pi(u)=x$. Since $u$ is a linear isomorphism of $\R^m$ onto
the tangent space $T_xM$, it induces an isomorphism of the tensor algebra over $\R^m$  onto the tensor
algebra over $T_xM$. Then $K_x$ is the image of $\KK$ under this isomorphism. The invariance of $\KK$ by $G$
implies that $K_x$ is defined independent of the choice of~$u$ in $\pi^{-1}(x)$. Then, we have the following
classical result
\cite{koba72}.
\begin{proposition}\label{prop:koba72}
Let $\KK$ be 
a tensor over the vector space $\R^m$
and $G$ the group of linear transformations 
of $\R^m$ leaving $\KK$ invariant\textup.
Let $P$ be a \Gdash structure on $M$ and $K$ the tensor 
field on $M$ defined by $\KK$ and $P$\textup. Then
\begin{enumerate}
\item[\textup{1.}] a diffeomorphism $\varphi \colon M\to M$ 
is an automorphism of the \Gdash structure~$P$ iff $\varphi$ leaves $K$ invariant\textup;
\item[\textup{2.}] a vector field~$\xi$ on~$M$ is an infinitesimal automorphism of~$P$ iff
$£_{L\xi}K=0$\textup.
\end{enumerate}
\end{proposition}
Now, we can use the concept of a \Gdash reductive Lie derivative to state an analogous result
for generalized automorphisms of~$P$.
\begin{proposition}\label{prop:god00}
In the same hypotheses of the previous proposition,
\begin{enumerate}
\item[\textup{1.}] an automorphism $\Phi \colon LM\to LM$ 
is a generalized automorphism of the \Gdash structure~$P$ iff $\Phi$ leaves $K$ invariant\textup;
\item[\textup{2.}] a \GLdash invariant vector field~$\Xi$ on~$LM$ is an 
infinitesimal generalized automorphism of~$P$ iff $£_\Xi K=0$\textup, whence $£_\Xi^GK\equiv 0$ for any 
\GLdash invariant vector field~$\Xi$ on~$LM$\textup.
\end{enumerate}
\end{proposition}
Note that the Lie derivative $£_{\Xi}^G K$ is well-defined since $K$ is a tensor field on~$M$
and therefore a section of a vector bundle associated with $W^{0,0}(LM)\cong LM$. Here, $Q=LM$ and
$H=\GL(m,\R)$. Nevertheless, consistently with what we said previously, $K$ has to be regarded
here as a section of a \emph{gauge-}natural, not simply natural, bundle over~$M$ (\cf\ Example~\ref{ex:btd}).
The choice $\Xi=L\xi$ reproduces Kobayashi's classical result.
\begin{corollary}
Let $\Xi$ be a \GLdash invariant vector field on~$LM$\textup, and let~$g$ be a metric tensor on~$M$ of
signature $(p,q)$\textup. Then\textup, $£_\Xi^{\SO(p,q)}g \equiv0$\textup.
\end{corollary}
The last corollary suggests that Bourguignon and Gauduchon's metric Lie derivative might be a
particular instance of a reductive Lie derivative. This is precisely the case, as explained in
the following fundamental
\begin{example}[Lie derivative of spinor fields. II] We know that the Kosmann lift~$\xiK$ onto
$\SO(M,g)$ of a vector field~$\xi$ on~$M$ is an \SOdash invariant vector field on $\SO(M,g)$,
and hence its lift~$\tilde\xiK$ onto $\Spin(M,g)$ is a \Spindash invariant vector field. As the spinor bundle
$S(M)$ is a vector bundle associated with $\Spin(M,g)$, the \SOdash reductive Lie derivative
$£_{L\xi}^{\SO(p,q)}\psi$ of a spinor field~$\psi$ coincides with $£_\xi\psi$, i.e.\ locally with
expression~\eqref{eq:dlK} or~\eqref{eq:dlKp}. Indeed, in this case we have, with an obvious notation,
$P=\SO(M,g)$, $G=\SO(M,g)$, $\tilde P=\Spin(M,g)$ and $\tilde P_{\tilde\l}=S(M)$.

\noindent
For $£_{L\xi}^{\SO(p,q)}g$ a similar remark to the one above for $£_{\Xi}^G K$ applies and
therefore, if $g=g_{\mu\nu}\,\d x^\mu\vee\d x^\nu$ in some natural chart, we have the local expression
\begin{align*}
£_{L\xi}^{\SO(p,q)}g_{\mu\nu} &\equiv\xi^\rho\de_\rho g_{\mu\nu}
+2g_{\rho(\mu}(\xiK)^\rho\ida_{\nu)}\\ &\equiv\xi^\rho\de_\rho g_{\mu\nu}
+g_{\rho(\mu}\de_{\nu)}\xi^\rho  -\delta^\rho\ida_{(\mu}g_{\nu)\sigma}\de_\rho\xi^\sigma
-\xi^\rho\delta^\sigma\ida_{(\mu|}\de_\rho g_{|\nu)\sigma}  \\ &\equiv0  \\
&\equiv £_\Xi^{\SO(p,q)}g_{\mu\nu},
\end{align*}
quite different from the usual (\emph{natural}) Lie derivative
\begin{align*}
£_{L\xi}g_{\mu\nu} &\equiv\xi^\rho\de_\rho g_{\mu\nu}
+2g_{\rho(\mu}(L\xi)^\rho\ida_{\nu)}  \\
&\equiv\xi^\rho\de_\rho g_{\mu\nu}
+2g_{\rho(\mu}\de_{\nu)}\xi^\rho  \\
&\equiv 2\cde_{(\mu}\xi_{\nu)}  \\
&\equiv£_\xi g_{\mu\nu}.
\end{align*}
\end{example}

\section*{Discussion}
In this paper we have investigated the hoary problem of the Lie derivative of spinor fields from a very
general point of view, following a functorial approach. We have done so by relying on three nice geometric
constructions: split structures, gauge-natural bundles and the general theory of Lie derivatives.

Such analysis has shown that, although for (purely) natural objects over a manifold~$M$ there is a
conceptually and mathematically natural definition of a Lie derivative with respect to a vector field
on~$M$, there is no such thing for more general gauge-natural objects, the vector field on~$M$ being
necessarily replaced by a \Gdash invariant vector field on some principal bundle $P(M,G)$. 

Conceptually speaking, though, there are at least two obvious definitions of a Lie derivative of spinor
fields, both relying on a canonical, not natural, lift of a vector field on~$M$ onto the bundle of its 
orthonormal frames, the so-called ``Kosmann lift''. Both definitions are geometrically well-defined and have
their own range of applicability, but, in general, only the gauge-natural one reduces to the standard
definition of a Lie derivative on natural objects.

\section*{Acknowledgements}
The authors are indebted to I.~Kol{\'a}{\v r} for many useful discussions and to J.~A. Vickers for a critical
reading of the manuscript. P.~M. acknowledges an EPSRC research studentship and a Faculty Research
Studentship from the University of Southampton.

\end{document}